\documentclass{article}

\usepackage{graphicx} 

\usepackage{amsmath,amssymb,amsthm}
\usepackage{mathrsfs}
\usepackage{tikz}
\usetikzlibrary{cd}
\usepackage{ upgreek }
\usepackage{ tipa }
\usetikzlibrary{cd}
\usepackage{dsfont}
\usepackage{extpfeil}
\usepackage[final]{todonotes}
\usepackage{xcolor}
\usepackage{bbold}
\usepackage{tikz}
\usetikzlibrary{cd}

\usepackage{hyperref}
\usepackage[notcite,notref,final]{showkeys}

\newcommand{\cC}{\mathcal{C}}

\newcommand{\cD}{\mathcal{D}}

\newcommand{\adjunction}[4]{\begin{tikzcd}[ampersand replacement=\&] #1 : #2 \arrow[r, shift left] \& #3 : #4 \arrow[l, shift left] \end{tikzcd}}

\newcommand{\module}[2][]{#2\text{-}\modop_{#1}}

\newcommand{\cA}{\mathscr{A}}
\newcommand{\chA}{\Ch(\cA)}

\DeclareMathOperator{\Ch}{Ch}
\DeclareMathOperator{\colim}{colim}
\DeclareMathOperator{\Ho}{Ho}
\DeclareMathOperator{\id}{id}
\DeclareMathOperator{\modop}{mod}
\DeclareMathOperator{\QCoh}{QCoh}

\newtheorem{theorem}{Theorem}[section]
\newtheorem{lemma}[theorem]{Lemma}

\newtheorem*{theorem*}{Theorem}
\newtheorem{proposition}[theorem]{Proposition}
\theoremstyle{definition}
\newtheorem{definition}[theorem]{Definition}

\newtheorem*{notation*}{Notations}

\theoremstyle{remark}
\newtheorem{remark}[theorem]{Remark}
\newtheorem{example}[theorem]{Example}
\newtheorem{corollary}[theorem]{Corollary}
\newtheorem{assumption}[theorem]{Assumptions}

\title{Stable compatibility for model categories}
\author{Kunal Arora \and Vivek Mohan Mallick}

\begin{document}

\maketitle
\begin{abstract}
By introducing the notion of stable compatibility, we provide a framework in the model category theory which addresses the existence of a stable model structure on the category of modules over a monad. Exact version of such monads have been studied in the setting of triangulated categories, separately as separable monads \`a la Balmer in tensor triangular geometry and non-separable monads \`a la Dell'Ambrogio and Sanders in stable homotopy theory. By introducing stable compatibility, we explore their properties together and produce a non-separable monad in Algebraic geometry and more.
\end{abstract}

\tableofcontents

\section{Introduction}
We are interested in the following broad question:
\begin{quote}
Given a stable model category  $\cC$ and a monad  $(M, \mu, \eta)$ acting on $\cC$, under what conditions is $\module{M}$, the category of modules of $M$ in $\cC$ i.e.\ the Eilenberg-Moore category, a stable model category?
\end{quote}
This interests us majorly in the case of Algebraic geometry where $\cC$ is $\Ch(QCoh(X))$, the category of chain complexes over quasi-coherent sheaves and $M$ is the monad generated by the adjunction induced by a map of quasi-compact and quasi-separated schemes $f:V \rightarrow X$.

The derived version of the above question is a result of Balmer (See, \cite{balmer:etalneth} theorem 3.5 for more details) about the monad $M=Rf_{*}f^{*}$ generated by a seperated \'etale morphism $f:V \rightarrow X$ of  quasi-compact and quasi-separated schemes acting on the unbounded derived category of quasi-coherent sheaves $D(QCoh(X))$ and it states that $\module{M}$  is triangulated and equivalent to $D(QCoh(V))$. 

This opens up the definition of \'etale morphism and hence \'etale-topology in the setting of \textit{tensor-triangular geometry}.

Our original motivation for this paper was to see what other morphisms and topologies from algebraic geometry  admit an extension in  \textit{tensor-triangular geometry}.

The monad associated to the seperated \'etale map is an example of \textit{separable monad}. We now state an example of a \textit{non-separable} monad from stable homotopy theory.

A similar situation arises in stable homotopy theory, when one considers the map say $f: R \rightarrow A$ between highly structured ring spectrum (say \cite{MR1417719} $S$-algebra), this \textit{change of ring} map induces an extension of scalars and restriction adjunction on the level of module categories $f^{*}=A \wedge - : M_{R} \leftrightarrows M_{A}:f_{*} $ where  $(f^{*} \dashv f_{*})$. Let $R$ be the sphere spectrum, denoted as $S$, and let $A$ be as above with $f: S \rightarrow A$ the unit map. The category of modules of the monad $f_{*}f^{*}$ on $M_{S}$ is $M_{A}$ which inherits a right-transferred model structure from $M_{S}$ and is simplicial, stable and cofibrantly generated model category.

Consider the monad formed by the derived adjunction $Lf^{*}: D(S) \leftrightarrows D(A): Rf_{*}$ between $D(S)=SH$ the stable homotopy of spectra and the derived category $D(A)$. Let $A\text{-Mod}_{SH}$ be the category of modules of the monad $Rf_{*}Lf^{*}$ on $D(S)$. Dell'Ambrogio and Sanders \cite[corollary 1.11]{MR3851535} characterizes the existence of triangulated structure on $A\text{-Mod}_{SH}$ in terms of its (exact) equivalence with $D(A)$. 
Example 1.12 (loc.\ cit.) says that Guti\'errez shows that corollay 1.11 holds for  $A=HB$ an Eilenberg Maclane spectrum of a field or subring B of $\mathbb{Q}$.

We note that the case of stable homotopy which provides an example of a non-separable monad says that we get a triangulation on  $A\text{-Mod}_{SH}$ only when it is equivalent to the homotopy category of a stable model category $M_{A}$. This lead us to investigate the question mentioned in the beginning and to provide a unified treatment for both separable and non-separable monads on stable model categories.

\subsection*{Main results}
Answering the above mentioned question packed in the definition of \textit{stable compatibility} (see definition \ref{def:stabcomp}), which arises when the forgetful functor preserves cofibrant objects under some suitable hypothesis, we have our first main result (see theorem \ref{the:stabcomp} below).
 
\begin{theorem*}
    
Let $\cC$ be a stable cofibrantly generated model category.
  Suppose $M$ is a monad on $\cC$ which satisfies the following.
  \begin{enumerate}
    \item
      The category $\module{M}$ admits a right transferred model structure, which is 
      cofibrantly generated, and such that the induced adjunction $\Sigma_{\module{M}}
      \dashv \Omega_{\module{M}}$ is a Quillen adjunction. 
    \item The monad $M$ satisfies stable compatibility.
  \end{enumerate}
  Then $\module{M}$ is a stable model category. Moreover, if the monad is induced by a
  Quillen adjunction $\adjunction{F}{\cC}{\cD}{G}$, the underlying comparison functor $K$
  is such that $K \circ F=F_{M}$ and $G=G_{M} \circ K$ as in the diagram below   
  \[
    \begin{tikzcd}
      \cC  \arrow[ddrr,"F_{M}", shift left=1mm] \arrow[dd,"F", shift left=1mm]
      &&
      \\
      && \\
      \cD \arrow[rr, dashrightarrow,"\exists !K" below] \arrow[uu,"G",shift left=1mm]
      && 
      \module{M} \arrow[uull, "G_{M}",shift left=1mm] 
    \end{tikzcd}
  \]
  and it preserves fibrations and acyclic fibrations.

\end{theorem*}
We use this result and its interaction with sections of a separable monad to obtain our second main result (see theorem \ref{the:mainsepm}).

\begin{theorem*}
    Let $\cC$ be a cofibrantly generated stable model category, and let $M$ be a 
separable monad on $\cC$ such that $M$ satisfies the assumptions of theorem 
\ref{the:stabcomp} via assumption \ref{asm:sepassmpt}. Then $\module{M}$ is a cofibrantly generated stable model 
category. Moreover, a cofiber sequence admits a direct summand inclusion into a suitable
cofiber sequence.
\end{theorem*}
Although this result provides a clarification to \cite[remark 6.1]{balmer:septrcat} but we mention it now that this does not address the question about the existence of an equivalence between $\Ho(\module{M})$ and $\module[\Ho \cC]{(\Ho M)}$ where the former is the homotopy category of $M$-modules and latter is the category of modules over the derived monad in the homotopy category .

This result provides us with a new model categorical insight about the structure of data of cofiber sequences (see remark \ref{rem:sepscmnd}), and for an application to the derived category $D(S)$ of a finite \'etale  $R$-algebra $S$ for a commutative ring $R$, see example \ref{exa:sepalgeb}.

Another strong application of \textit{stable compatibility} is in algebraic geometry where we show that finite and flat morphisms between suitable schemes form a monad whose module category turns out to be a cofibrantly generated stable model category (see \textsection \ref{ssc:finflmor}). In this case, we even have an equivalence between between $\Ho(\module{M})$ and $\module[\Ho \cC]{(\Ho M)}$. This is a first major step outside \'etale morphisms in \textit{tensor triangular geometry}.

We mention other potential applications and conclude introduction in this section-wise overview. 

Section \ref{sec:prereqst} provides the preliminaries on Eilenberg-Moore categories $\module{M}$, simplicial model categories, the pre-triangulated structure in the sense of \cite{MR1650134} and stable model categories. We also mention about the existence of right-transferred model structure on $\module{M}$, and lastly we recall the locally free model structure on $\Ch(QCoh(S))$ of a suitable scheme. 

Section \ref{sec:stablecpt} is our main section on \textit{Stable compatibility} with our main result theorem \ref{the:stabcomp} and an auxiliary lemma \ref{lem:cndforsc}, which helps us in detecting scenarios where \textit{stable compatibility} holds. It forms our main tool in \textsection \ref{ssc:finflmor}. An analogous result, concluding stability of the module category in the $\infty$-category setup, can be deduced from \cite[corollary 5.5 and theorem 5.7]{MR3338539}. 

Section \ref{sec:sepmonad} addresses the interaction of \textit{stable compatibilty} with the definition of a \textit{separable monad}.

Section \ref{sec:applications} mentions other potential applications in stable homotopy (see \textsection \ref{ssc:drumhack}), which comes along with an improved observation to a result of [\cite{MR3937590} Theorem 2.3]. Lastly, we mention that in the dual setup of \textit{stable compatibility}, one can recover the stability of the category of comodules over a comonad in \textsection \ref{ssc:comodules}.

\subsection*{Acknowledgements} Both the authors thank IISER Pune for funding and providing
infrastructure to conduct this research. This research was supported in part by the International Centre for Theoretical Sciences (ICTS) for participating in the program -  Dualities in Topology and Algebra (code: ICTS/dta2023/05). 

\section{Prerequisites} \label{sec:prereqst}

\subsection{Module categories} \label{ssc:monads}
Throughout this paper, a monad on a category $\cC$ is denoted by a triple $(M,\mu,\eta)$ (or simply $M$, when the context is clear). By $\module[\cC]{M}$ (or simply, $\module{M}$ when there is no ambiguity), we denote the category of $M$-modules, the Eilenberg-Moore category \cite{maclane:catworkmath} associated with $M$. The free-forgetful adjunction is denoted by 
\[
\adjunction{F_M}{\cC}{\module[\cC]{M}}{G_M}.
\]
The following definitions are taken from  \cite{balmer:septrcat}.

\begin{definition}(\cite{balmer:septrcat}, Definition 2.1)\label{def:monad}  Let $\cC$ be a category. A \textit{monad} on $\cC$ is a triple $(M,\mu,\eta)$ where $M: \cC \rightarrow \cC$ is an endofunctor, $\mu:M^{2} \rightarrow M$ (\textit{the multiplication}) and $\eta:Id_{\cC} \rightarrow M$ (\textit{the unit}) are natural transformations such that the following diagrams expressing associativity and two-sided unit,commute:
\begin{equation}
    \begin{tikzcd}
      M^{3} \arrow{r}{M\mu} \arrow[swap]{d}{\mu M} & M^{2} \arrow{d}{\mu} \\
      M^{2} \arrow{r}{\mu} & M
      \end{tikzcd}
     and 
     \begin{tikzcd}
      M \arrow{r}{M\eta} \arrow[swap]{dr}{id} & M^{2} \arrow{d}{\mu} & M \arrow[swap]{l}{\eta M} \arrow{dl}{id} \\
       &M .
     \end{tikzcd}
\end{equation}
When the category $\cC$ is suspended(i.e. an additive category with an autoequivalence $\Sigma:\cC \xrightarrow{\sim}\cC$, see \cite{balmer:septrcat}  Definition 1.1 for more), we say that an additive monad $M$ is \textit{stable} if $M$, $\mu$, $\eta$ commute with suspension: $\Sigma M=M \Sigma , \mu \Sigma=\Sigma \mu , \eta\Sigma=\Sigma\eta.$
\end{definition}
See also the assumptions \ref{asm:stabcomp} below.
\begin{example}
    (See, \cite{balmer:septrcat} Example 2.3) Let $\cC$ be a monoidal category with tensor $\bigotimes: \cC \times \cC \rightarrow \cC$ and the unit $\mathbb{1}  \in \cC$. A (unital and associative) ring object in $\cC$ is a triple $(A,\mu,\eta)$ where $A$ is an object of $\cC$ and where the \textit{multiplication} $\mu: A \bigotimes A \rightarrow A$ and the unit $\eta: \mathbb{1} \rightarrow A$ are morphisms in $\cC$ satisfying the usual associativity and two-sided unit conditions analogous to Definition 2.1 (1). Let $M:\cC \rightarrow \cC$ be the functor $A \bigotimes-$, with the obvious $\mu: M^{2} \rightarrow M$ and $\eta:Id \rightarrow M.$ Then $M$ is a monad on $\cC$.

\end{example}
For a commutative ring $R$ and for $\cC=R-Mod$ with $\bigotimes=\bigotimes_{R}$, the ring objects in $\cC$ are the usual $R$-algebras. In particular, ring objects in $\mathbb{Z}-Mod$ are ordinary rings. So, monads are generalizations of rings and algebras.

\begin{definition}(See, \cite{balmer:septrcat} Definition 2.4)\label{def:M-mod} Let $M:\cC \rightarrow \cC$ be a monad. The \textit{Eilenberg-Moore category} of (left)$M$-modules $M-mod$ is defined as follows. A left $M$-module is a pair $(x, \lambda)$ where $x$ is an object of $\cC$ (\textit{the underlying object}) and $\lambda:Mx \rightarrow x$ is a morphism(\textit{the left action}) such that the following diagrams both commute:
\[\begin{tikzcd}
M^{2}x \arrow{r}{M\lambda} \arrow[swap]{d}{\mu_{x}} & Mx \arrow{d}{\lambda} \\
Mx \arrow{r}{\lambda} & x
\end{tikzcd}
\] and 
\[\begin{tikzcd}
x \arrow{r}{\eta_{x}} \arrow[swap]{dr}{id} & Mx \arrow{d}{\lambda} \\
& x
\end{tikzcd}
\]
A morphism $f: (x,\lambda) \rightarrow (x^{'}, \lambda')$ of left $M$-modules is a morphism $f: x \rightarrow x^{'}$ in $\cC$ which is $M$-linear, i.e. such that $\lambda'\circ M(f)=f \circ \lambda$. The \textit{free module} functor $F_{M}: \cC \rightarrow M-Mod$ defined by $F_{M}(y):=(My,\mu_{y})$. The right adjoint $G_{M}:M-Mod \rightarrow \cC$ just forgets the action. We define the \textit{Klesli category of free $M$-modules} $M-Free_{\mathcal{C}}$ as the full subcategory $F_{M}(\cC)$ of  $M-Mod$. The above functors $F_{M}$ and $G_{M}$ restrict to an adjunction between $\cC$ and $M-Free_{\cC}$: \[\begin{tikzcd}
 & \cC \arrow[dl,"F_{M}"] \arrow[dr,"F_{M}"]  \\
 M-Free_{\cC} \arrow[hookrightarrow, rr,"(fully faithful)"below] \arrow[ur,"G_{M}",shift left=1ex] && M-Mod \arrow[ul, "G_{M}",shift left=1ex] 
\end{tikzcd}
\]
When $\mathcal{C}$ is suspended and $M$ is stable, then both $M-Free_{\mathcal{C}}$ and $M-Mod$ inherit an obvious suspension such that $F_{M}$ and $G_{M}$ commute with that suspension.

\end{definition}
Contrast the \textit{stable monad} with the  assumptions \ref{asm:stabcomp} and Remark 3.9 below.
\begin{example}
    (See, \cite{balmer:septrcat} Example 2.5) When $(A,\mu,\eta)$ is a ring object in a monoidal category $\cC=(\cC, \bigotimes,\mathbb{1})$ as in example 2.2, the above constructions yield the natural categories of left $A$-modules and free left $A$-modules.
\end{example}
\begin{remark}(See, \cite{balmer:septrcat} Remark 2.7)
    Given an adjunction $F:\mathcal{C} \leftrightarrows \mathcal{D}: G$, let $\eta: Id_{\mathcal{C}} \rightarrow GF$ be the unit and $\epsilon:FG \rightarrow Id_{\mathcal{D}}$ be the counit of this adjunction. Then $M:=GF$ is a monad on $\mathcal{C}$, with unit $\eta$ and multiplication $\mu:=G\epsilon F:M^{2}=G(FG)F \xrightarrow{G \epsilon F} GF=M$. One says that $M$ is \textit{realized by the adjunction} $(F,G)$.

    If moreover $\mathcal{C}$ and $\mathcal{D}$ are suspended and $F$ and $G$ commute with suspension, we tacitly assume that the adjunction respects the suspension, meaning that $\eta$ and $\mu$ commute with suspension. In that case, the monad $GF$ is stable.
\end{remark}

Given a monad $M:\cC \rightarrow \cC$ , there are in general many adjunctions realizing $M$. They form a category in which the Klesli construction $M-Free_{\mathcal{C}}$ is initial and the Eilenberg-Moore category $M-Mod$ is final, as stated precisely below.

\begin{proposition}[\cite{balmer:septrcat}, proposition 2.8] \label{pro:kleisliem}
Let $M:\cC \rightarrow \cC$  be a monad  \textit{realized by an adjunction }$F:\mathcal{C} \leftrightarrows \mathcal{D}: G$. Then there are unique functors L and K as follows: \[\begin{tikzcd}
 & \cC \arrow[dl,"F_{M}"] \arrow[dr,"F_{M}"] \arrow[d,"F"] \\
 M-Free_{\cC} \arrow[dashrightarrow,r,"\exists !L"below] \arrow[ur,"G_{M}",shift left=1ex] & D \arrow[dashrightarrow,r,"\exists !K"below] \arrow[u,"G",shift left=1ex]& M-Mod \arrow[ul, "G_{M}",shift left=1ex] 
\end{tikzcd}
\]
such that $L \circ F_{M}=F, G_{M}=G \circ L$
and $G=G_{M} \circ K$. Moreover, $L$ is fully faithful. 
Finally, if $\mathcal{C}$ and $\mathcal{D}$ are suspended, if $M$ is stable and if $(F,G)$ is an adjunction of functors commuting with suspension, then the functors $L$ and $K$ commute with suspension as well.

\end{proposition}
Contrast this with the Theorem \ref{the:stabcomp}below.
\subsection{Model categories}
Throughout this paper our model categories are as defined in \cite{MR1650134} definition 1.1.3 (See for more details). The exciting feature of this definition is that the functorial factorizations are a part of the model structure. 
\begin{definition}\label{def:pointedmodcat}
    (As on page 4, \cite{MR1650134}) We call a model category(or any category with an initial and terminal object) \textit{pointed} if the map from the initial object to terminal object is an isomorphism. Given a model category $\cC$, define $\cC_{*}$ to be the category under the terminal object $*$. That is, an object of $\cC_{*}$ is a map $* \xrightarrow{v} X$ of $\cC$, often written $(X,v)$. We think of $(X,v)$ as an object together with a \textit{basepoint} $v$. A morphism from $(X,v)$ to $(Y,w)$ is a morphism $X \rightarrow Y$ of $\cC$ that takes $v$ to $w$.
\end{definition}
By (\cite{MR1650134}, Proposition 1.1.8) one can induce a model structure on $\cC_{*}$ from $\cC$ via the adjunction forgetful functor $U: C_{*} \rightarrow C$ and it's left adjoint of adding the basepoint i.e. X goes to $X_{+}=X \sqcup *$.
\\Let's recall some definitions needed to define \textit{Cofibrantly generated model categories}
\begin{definition}
    (\cite{MR1650134}, Definition 2.1.3) Suppose $\cC$ is a category with all small colimits, $\cD$ is a collection of morphisms in $\cC$, $A$ is an object of $\cC$ and $\kappa$ is a cardinal. We say that $A$ is $\kappa$-\textit{small relative} to $\cD$ if, for all $\kappa$-filtered ordinals $\lambda$ and all $\lambda-$sequences $X_{0} \rightarrow X_{1} \rightarrow...\rightarrow X_{\beta}\rightarrow ...$ such that each map $X_{\beta} \rightarrow X_{\beta+1}$ is in $\cD$ for $\beta+1 < \lambda$, the map of sets $\colim_{\beta < \lambda} \cC(A,X_{\beta}) \rightarrow \cC(A, \colim_{\beta < \lambda} X_{\beta})$ is an isomorphism. We say that $A$ is \textit{small relative} to $\cD$ if  it is $\kappa-$small relative to $\cD$ for some $\kappa$. We say that $A$ is small if it is small relative to $\cC$ itself.
\end{definition}
From (\cite{MR1650134}, Definition 2.1.7) we state that for $I$ a class of maps in $\cC$, a map is \textit{$I$-injective} if it has the right lifting property with respect to every map in $I$. The class of \textit{$I$-injective} maps is denoted by $I$-inj. Similarly, a map with the left lifting property with respect to the class $I$-inj is called \textit{$I$-cofibration} and the class of such maps is denoted by $I$-cof. 
\begin{definition}
    (\cite{MR1650134}, Definition 2.1.9)\label{def:I-cell} Let $I$ be a set of maps in a category $\cC$ containing all small colimits. A \textit{relative $I$-cell complex} is a transifinite composition of pushouts of elements of $I$ i.e. if $f: A \rightarrow B$ is a relative $I$-cell complex, then there is an ordinal $\lambda$ and a $\lambda$-sequence, $X \colon \lambda \longrightarrow \cC$, such that $f$ is the composition of $X$ and such that, for each $\beta$ with $\beta+1 < \lambda$, there is a pushout square
     \[\begin{tikzcd}
        C_{\beta} \arrow{r} \arrow{d}{g_{\beta}}  & X_{\beta} \arrow{d} \\
        D_{\beta} \arrow{r} & X_{\beta+1}
     \end{tikzcd}
     \]  
     such that $g_{\beta} \in I$. We denote the collection of relative $I$-cell complexes by \textit{$I$-cell}.
\end{definition}
We now state (\cite{MR1650134}, Theorem 2.1.14) in the form of a lemma.
\begin{lemma}
    (The small object argument) Suppose $\cC$ is a category containing all small colimits, and $I$ is a set of maps in $\cC$. Suppose the domains of the maps of $I$ are small relative to $I$-cell. Then there is a functorial factorization $(\gamma,\delta)$ on $\cC$ such that, for all morphisms $f$ in $\cC$, the map $\gamma (f)$ is in $I$-cell and the map $\delta(f)$ is in $I-inj$. 
\end{lemma}
We now recall the definition of a \textit{Cofibrantly generated model category}. 
\begin{definition} \label{def:cofgenemodcat}
    (\cite{MR1650134}, Definition 2.1.17) Suppose $\cC$ is a model category. We say that $\cC$ is \textit{cofibrantly generated} if there are sets $I$ and $J$ of maps such that:
    \begin{enumerate}
        \item The domains of the maps of $I$ are small relative to $I$-cell;
        \item The domains of the maps of $J$ are small relative to $J$-cell;
        \item The class of fibrations is $J-inj$; and 
        \item The class of trivial fibrations is $I-inj$.
    \end{enumerate}
    We refer to $I$ as the set of \textit{generating cofibrations}, and to $J$ as the set of \textit{generating trivial cofibrations}.
\end{definition}
We note that the small object argument of lemma 2.10 gives the functorial factorizations in a \textit{cofibrantly generated model category}. For more details see (\cite{MR1650134}, Proposition 2.1.18) and also Theorem 2.1.19).

\subsection*{Simplicial model categories}

We now recall the notion of a \emph{pointed simplicial model category}. Let $SSet$ be the category of simplicial sets and $SSet_{*}$ be the pointed simplicial sets.
\begin{definition}
    (See \cite{MR1650134}, Definition 4.2.18) Given a monoidal model category $\cC$, a $\cC$-\textit{model category}  is a $\cC$-module $\cD$(see Definition 4.1.6 for more details) with a model structure making $\cD$ into a model category such that the following conditions hold.
    \begin{enumerate}
        \item The action map $\bigotimes:\cC \times \cD \rightarrow \cD$ is a Quillen bifunctor.
        \item If $QS \xrightarrow{q} S$ is the cofibrant replacement for $S$ in $\cC$, then the map $X \bigotimes QS  \xrightarrow{1\bigotimes q} X \bigotimes S$ is a weak equivalence for all cofibrant $X$ $\in \cD$.
    \end{enumerate}
\end{definition}
\begin{definition}\label{def:pointedsimpmodcat}
    An $SSet$-model category is called a \textit{simplicial model category}. a pointed simplicial model category is the same thing as $SSet_{*}-$model category.
\end{definition}
We recall a Remark which concludes that the simplicial suspension and loop functor( i.e. the pointed tensor and cotensor with simplical circle $S^{1}=\triangle[1]_{+}/ \partial \triangle[1]_{+1})$ is a Quillen pair.
\begin{remark}
    (See \cite{MR1650134}, Remark 4.2.3) Suppose $\bigotimes:\cC \times \cD \rightarrow \mathcal{E}$ is a Quillen bifunctor. Then, if C is cofibrant, the functor $C \bigotimes -: \cD \rightarrow \mathcal{E}$ is a Quillen functor with right adjoint $Hom_{l}(C,-)$. Similarly if D is cofibrant , the functor $- \bigotimes D$ is a Quillen functor with right adjoint $Hom_{r}(D,-)$. Dually we have a statement with fibrant objects which we don't need.
\end{remark}
\begin{corollary}
    
Since a pointed simplical circle is a cofibrant object in $SSet_{*}$, so we get that the simplical suspension and loop functor is a Quillen pair on a pointed simplicial model category.
\end{corollary}
By the techniques of $framing$ in (See \cite{MR1650134}, Ch-5) Hovey shows that the homotopy category of a model category is naturally a closed $Ho SSet$-module and the homotopy category of a pointed model category is naturally a closed $Ho-SSet_{*}$-module.
And as noted in section 5.7 corollary 5.7.7 that a pointed simplicial model structure induces this closed $Ho-SSet_{*}$-module structure on the homotopy category of a pointed simplicial model category via the homotopy pseudo-2-functor.
\\This lets us give the following definition for homotopy category of any pointed model category. 
\begin{definition}
    (See \cite{MR1650134}) Definition 6.1.1) Suppose $\cC$ is a pointed model category. The \textit{suspension functor} $\Sigma: Ho \cC \rightarrow Ho\cC$ is the functor $X \mapsto X \wedge^{L} S^1$ defined by the closed action of $HoSSet_{*}$ on Ho$\cC$ as given in section 5.7. Dually the loop functor $\Omega: Ho \cC \rightarrow Ho\cC$ is the functor $X \mapsto RHom_{*}(S^{1},X)$.
\end{definition}
\begin{corollary} \label{cor:Quiequvforsimpstab}
    For a pointed simplicial model category the definition 2.16 is the one given by deriving the Quillen adjunction of corollary 2.15.
\end{corollary}

\subsection*{Pre-triangulated structure on the homotopy category of a model category}
We now recall the definition of cofiber and fiber sequences.

\begin{definition}
    (See \cite{MR1650134}, Definition 6.2.6)\label{def:cofseq} Suppose $\cC$ is a pointed model category. A \textit{cofiber sequence} in Ho $\cC$ is a diagram $X \xrightarrow{f} Y \xrightarrow{g} Z$ in Ho$\cC$ together with a right coaction of $\Sigma X$ on $Z$ which is isomorphic in Ho$\cC$ to a diagram of the form $A \rightarrow B \rightarrow C$ where $f$ is a cofibration of cofibrant objects in $\cC$ with cofiber $g$ and where $C$ has the right $\Sigma A$-coaction(This will give a coboundary map as in defintion 2.19 below) given by theorem 6.2.1. Dually a \textit{fiber sequence} is a diagram $X \rightarrow Y \rightarrow Z$ together with a right action of $\Omega Z$ on $X$ which is isomorphic to a diagram $F \xrightarrow{i} E \xrightarrow{p} B$ where $p$ is a fibration of fibrant objects with fiber $i$ and where $F$ has the right $\Omega B$-action given by theorem 6.2.1.
\end{definition}
\begin{definition} \label{def:boundarymap}
    (See \cite{MR1650134}, Definition 6.2.7) Suppose $\cC$ is a pointed model category, and $X \xrightarrow{f} Y \xrightarrow{g} Z$ is a cofiber sequence in Ho$\cC$. The \textit{boundary map} is the map $\partial: Z \rightarrow \Sigma X$ in Ho$\cC$ which is the composite $Z \rightarrow Z\bigsqcup \Sigma X \xrightarrow{0 \times 1} \Sigma X$ \\ where the first map is the coaction as mentioned above. Dually, if $X \xrightarrow{f} Y \xrightarrow{g} Z$ is a fiber sequence, the \textit{boundary map} is the map $\partial :\Omega Z \rightarrow X$ which is the composite $\Omega Z \xrightarrow{(0,1)} X \times \Omega Z \xrightarrow{0 \times 1}  X$
\end{definition}
We now define what is called a pre-triangulated structure on a closed $HoSSet_{*}$-module, hence on the homotopy category of a pointed model category.
\begin{definition}
    (See \cite{MR1650134}, Definition 6.5.1) \label{def:preTristruct}Suppose $\mathcal{S}$ is a nontrival (right) closed $Ho SSet_{*}$-module. A \textit{pre-triangulation} on $\mathcal{S}$ is a collection of \textit{cofiber sequences}, or \textit{left triangles}, and \textit{fiber sequences}, or \textit{right triangles}, satisfying the following conditions.
    \begin{enumerate}

        \item A cofiber sequence is in particular a diagram of the form $X \xrightarrow{f} Y \xrightarrow{g} Z$ together with a right coaction of the cogroup $\Sigma X$ on $Z$. A fiber sequence is in particular a diagram of the form $X \xrightarrow{f} Y \xrightarrow{g} Z$ together with a right action of the group $\Omega Z$ on $X$.
        \item Every diagram isomorphic to a cofiber sequence is a cofiber sequence. Similarly, every diagram isomorphic to a fiber sequence is a fiber sequence. The isomorphisms respect the action as in \cite{MR1650134}, lemma 6.3.1.
        \item For any $X$, the diagram $* \rightarrow X \xrightarrow{1} X$ is a cofiber sequence(with the only possible coaction). The diagram $X \xrightarrow{1} X \rightarrow *$ is a fiber sequence (with the only possible coaction).
        \item Every map is a part of a cofiber fiber sequence, as in \cite{MR1650134}, Lemma 6.3.3.
        \item Cofiber sequences can be shifted to the right, and fiber sequences can be shifted to the left, as in \cite{MR1650134}, Proposition 6.3.4.
        \item Fill-in maps exist, as in \cite{MR1650134}, Proposition 6.3.5.
        \item Verdier's octahedral axiom and its dual both hold, as in \cite{MR1650134}, Proposition 6.3.6.
        \item Cofiber and fiber sequences are comptaible, as in \cite{MR1650134}, Proposition 6.3.7.
        \item The smash product preserves cofiber sequences in each variable. The functor $Hom(-,-)$ preserves fiber sequences in the second variable and converts cofiber sequences in the first variable into fiber sequences. The functor $Map(-,-)$ preserves fiber sequences in the second variable and converts cofiber sequences in the first variable into fiber sequences. See \cite{MR1650134}, corollary 6.4.2 and \cite{MR1650134}, Proposition 6.4.5 for exact statements.
    \end{enumerate}
\end{definition}
\begin{definition}
    (See, \cite{MR1650134}, Definition 6.5.2)\label{def:pretriCat} A \textit{pre-triangulated category} is a nontrivial closed $HoSSet_{*}$ modules $\mathcal{S}$ with all small coproducts and products, together with a pre-triangulation on $\mathcal{S}.$
\end{definition}
\begin{remark}
    As we note this \textit{pre-triangulated structure} exists on the homotopy category of a model category and a given Quillen adjunction say $(F \dashv G)$ induces an \textit{exact} adjunction $(LF \dashv RG)$between homotopy categories, where $LF$ preserves cofiber sequences and $RG$ preserves fiber sequences as in definition (\cite{MR1650134}, Definition 6.5.4). 
\end{remark}
 We now define a Stable model category.
 \begin{definition}\label{def:Hoveytricat}
     (See, \cite{MR1650134}, Definition 7.1.1) A \textit{triangulated category} is a pre-triangulated category in which the suspension functor $\Sigma$ is an equivalence of categories.\label{def:stabmodcat} A pointed model category is \textit{stable} if it's homotopy category is triangulated.
 \end{definition}
 \begin{example}
     All the categories in this paper such as $Ch(\mathcal{A})$, where $\mathcal{A}$ is any Grothendieck abelian category, \cite{MR1417719} $S$-modules, etc. are stable.
 \end{example}
 \begin{remark} \label{rem:factstri}
 \begin{enumerate}
     \item \label{ite:stabquilleneq} A simplicial model category is stable if it is pointed and the Quillen adjunction of corollary 2.15 is a Quillen equivalence.
     \item \label{ite:hoveyvsverdier} By proposition (\cite{MR1650134}, 7.1.6) a triangulated category in the sense of definition 2.23 with \textit{suspension functor} and \textit{cofiber sequences} becomes a classical triangulated category in the sense of \cite{neeman:trcat}.
     \item Although \cite{MR1417719} $S$-modules is also an example of simplicial stable cofibrantly generated model category. In folklore $Ch(\mathcal{A})$ can also be made simplicial. But we note that we don't need simplicial for our main result Theorem 3.5 and hence not for $Ch(\mathcal{A})$ as well.
\end{enumerate}
 \end{remark}

\subsection{Transferred model structure on module categories.}
We now state a result, which makes sure $\module{M}$ admits a \textit{right-induced} model structure from $\cC$.
Define a map of $\module{M}$ to be a \textit{weak equivalence} or \textit{fibration} if the underlying map in $\cC$ is a \textit{weak equivalence} or \textit{fibration} respectively. Define a map in $M$-mod to be a \textit{cofibration} if it has the left lifting property with respect to all acyclic fibrations.
A \textit{regular $I$-cofibration} is a map in the class of the (possibly transfinite) compositions of pushouts of maps in $I$.

Let $X$ be a $M$-module. Define a \textit{path object} for $X$ to be an object in $M$-mod $X^{I}$ together with $M$-algebra maps $X \xrightarrow{\sim} X^{I} \rightarrow X \times X$ factoring the diagonal map, such that the first map is a weak equivalence and the second map is a fibration in the underlying category $\cC$.

The following is lemma 2.3 from \cite{MR1734325} (\emph{Caveat:} In that paper, $M$-algebras are what we call $M$-modules here).
\begin{lemma} \label{lem:schwship}
    Let $\mathcal{C}$ be a cofibrantly generated model category, with $I$ a set of generating cofibrations and $J$ a set of generating acyclic cofibrations. Let $(M,\mu,\eta)$ be a monad on $\mathcal{C}$ such that $M:\mathcal{C} \rightarrow \mathcal{C}$ commutes with filtered direct limits. Let $I_{M}$ and $J_{M}$ denote the images of $I$ and $J$ under the free $M$-module functor $F_{M}$. Assume that the domains and codomains of the morphisms in $I_{M}$ and $J_{M}$ are small with respect to regular $I_{M}$-cofibrations and regular $J_{M}$-cofibrations respectively. Suppose that  
    \begin{enumerate}
        \item every regular $J_{M}$-cofibration is a weak equivalence, or
        \item every object of $\cC$ is fibrant and every $M$-module admits a path object.
    \end{enumerate}
    Then the category M-mod is a cofibrantly generated model category with $I_{M}$ a generating  set of cofibrations and $J_{M}$ a generating set of acyclic cofibrations,
\end{lemma}
 The forgetful functor $G_{M}:M-mod \rightarrow \cC$ is said to \textit{right-induce} a cofibrantly generated model category structure on $M-mod$.
 \begin{remark} \label{rem:postcofib}
     \begin{enumerate}
         \item As mentioned in \cite{MR1734325} remark 2.4, when a category $\cC$ is a \textit{locally presentable category} every object of it is small. Further, when the category is not locally presentable then we can rely on the assumption that if underlying functor $M$ commutes with filtered direct limits, then so does $G_{M}$. Hence by adjointness, if every object of $\cC$ is small relative to $\cC$, then every free $M$-module is small relative to $M-mod$, so the smallness assumptions of lemma 2.26 holds. 
         \item \label{rem:simppath} As noted in \cite{MR1734325} remark 4.5, if $\cC$ is a simplical model category on which a simplical monad acts then M-mod obtains a simplicial structure in which case it always has path objects. In cases such as \cite{MR1417719} $S$-modules all objects are fibrant as well.
         \item Another set of conditions for which lemma 2.26 holds is presented in lemma \ref{lem:cndforsc} below.
     \end{enumerate}
 \end{remark}

\subsection{Algebraic Geometry}
\begin{notation*}
A Grothendieck abelian category $\mathcal{A}$, we consider in this paper, is the category of \textit{quasi-coherent sheaves} over a scheme $S$ denoted by $QCoh(S)$. We denote by $\Ch(QCoh(S))$ the unbounded chain complexes of quasi-coherent sheaves and it's homotopy category the derived category is denoted by $D(QCoh(S))$.

We now recall a model structure in Algebraic geometry called \textit{locally free model structure} as presented by Hovey \cite{MR1814077} on $\Ch(QCoh(S))$ for a suitable scheme $S$ as below. \\ 
Let $\mathcal{A}$ be a Grothendieck category and recall that every object in a Grothendieck category is small. This lets us deal with smalless in transferred model structures as in Remark 2.27 (1) above.

Recall that $S^{n-1},D^{n}:\mathcal{A} \rightarrow Ch(\mathcal{A})$ are the functors that send $A \mapsto S^{n-1}A$ i.e. the complex whose only non-zero object is $A$ in degree $n-1$ and $A \mapsto D^{n}A$ the complex whose only nonzero objects are $A$ in degree $n-1$ and $n$ respectively.
\end{notation*}

Recall the following definition from \cite{MR1814077}.
\begin{definition}
    (\cite{MR1814077}, Definition 1.1) \label{def: I and J}Suppose $\mathcal{M} $ is a set of monomorphisms in a Grothendieck category $\mathcal{A}$. Let $\mathcal{F}$ denote the set of codomains of the maps of $\mathcal{M}.$ We will say that $\mathcal{M}$ is pointed if $0 \in \mathcal{F}$ and, if $F \in \mathcal{F}$, then $0 \rightarrow F$ is in $\mathcal{M}.$ Define $J$ to be the set of all $D^{n}f$, where $n$ is an integer and $f \in \mathcal{M}$. Define $I$ to be the union of $J$ and the maps $S^{n-1}F \rightarrow D^{n}F$ for $F \in \mathcal{F}$ and $n$ an integer. Then define a map $p$ to be an $\mathcal{M}\text{-fibration}$ if $p$ is in $J\text{-inj}$, and define $p$ to be an $\mathcal{M}$-cofibration if $p$ is in $I\text{-cof}$. 
\end{definition}

Recall that an object $B$ of a Grothendieck category $\mathcal{A}$ is said to have \textit{finite projective dimension} if there is an integer $n_{0}$ such that $Ext_{\mathcal{A}}(B,C)=0$ for all $n \ge n_{0}$ and all objects $C$ of $\mathcal{A}$. 

We now recall a result about existence of a cofibrantly generated model structure on $Ch(\mathcal{A})$ .
\begin{theorem} (\cite{MR1814077} Theorem 2.2)
    Suppose $\mathcal{A}$ is a Grothendieck category with a set $\mathcal{F}$ of generators, each of which has finite projective dimension. Let $\mathcal{M}$ denote the set of inclusions $A \rightarrow F$ of subobjects of objects $F \in \mathcal{F}$. Then there is a proper cofibrantly generated model structure on $Ch(\mathcal{A})$, where the weak equivalences are the quasi-isomorphims, the fibrations are the dimensionwise split surjections with dimensionwise injective kernel, and the cofibrations are the $\mathcal{M}$-cofibrations.
\end{theorem}

Also recall that a Noetherian scheme $S$ is said to have \textit{enough locally frees} if every coherent sheaf on $S$ is a quotient of a locally free sheaf of finite rank. 

We now recall  \cite[proposition 2.3]{MR1814077} which allows the existence of the mentioned \textit{locally free model structure}.

\begin{proposition}[\cite{MR1814077}, Proposition 2.3] Suppose S is a noetherian scheme with enough locally frees. In addition, suppose that S is either finite-dimensional or is separated. Then the locally free sheaves of finite rank are generators of finite projective dimension for the category QCoh(S) of quasi-coherent sheaves on S.
\end{proposition}

We have the following existence theorem as a corollary to Thereom 2.29 and proposition 2.30, which is \cite{MR1814077} Theorem 2.4.
\begin{theorem} (\cite{MR1814077} Theorem 2.4)\label{thm:locallyfreeModelStruct}
    Suppose $S$ is a noetherian scheme with enough locally frees, and suppose that $S$ is either finite-dimensional or seperated. Then there is a proper,cofibrantly generated, model structure on the category $\Ch(QCoh(S))$ of unbounded complexes of quasi-coherent sheaves, where the weak equivalences are the quasi-isomorphisms and the fibrations are the dimensionwise split surjections with dimensionwise injective kernel.
\end{theorem}

\begin{remark}
\begin{enumerate}
    \item If $f:X \rightarrow Y$ is a map  between schemes satisfying the hypotheses of thereom 2.31 then $f^{*}: QCoh(Y) \rightarrow QCoh(X)$ will induce a Quillen functor between the locally free model structures if and only if $f^{*}$ is exact. Note that for converse $f^{*}$ preserves locally free sheaves of finite rank.
    
    \item We have the following two deductions about $D(QCoh(S))$ from the locally free model structure: 
    \begin{enumerate}
        \item \cite[corollary 2.5]{MR1814077} Suppose $S$ is a noetherian scheme with enough locally frees, and $S$ is either finite-dimensional or seperated. Then the locally free sheaves of finite rank form a set of small weak generators for the $D(QCoh(S))$.
        \item \cite[corollary 2.6]{MR1814077} Suppose $S$ is a finite-dimensional noetherian scheme with enough locally frees. Then the category $D(QCoh(S))$ is a unital algebraic stable homotopy category, where the generators are the locally free sheaves of finite rank.             
    \end{enumerate}  
\end{enumerate}

\end{remark}

\section{Stable Compatibility} \label{sec:stablecpt}
In this paper, we assume the following:
\begin{assumption} \label{asm:stabcomp}
    \begin{enumerate}
        \item The category $\cC$ is pointed.
        \item We work with monads $(M, \mu, \eta)$ (see \ref{ssc:monads}), where the underlying endofunctor $M: \cC \rightarrow \cC$ of the monad preserves the zero object, and hence, the category of $M$-modules $\module[\cC]{M}$ is also pointed. The free functor $F_M$ and its right adjoint $G_M$ map the zero object in the appropriate category to the zero object
        of the other.
        \item \label{itm:suslooqe} Moreover, $\cC$ is a cofibrantly generated stable model category such that $(\Sigma \dashv \Omega)$ is a Quillen equivalence. This is the case, for example, when the cofibrantly generated stable model category is simplicial as noted in remark \ref{rem:factstri}\eqref{ite:stabquilleneq}.
        \item \label{itm:suscomGM} We also assume that $M$ commutes with the suspension $\Sigma$.
        \item \label{itm:rtmodstr} Moreover, when $\module{M}$ admits a right-transferred model structure, we shall assume that 
        $\adjunction{F_M}{\cC}{\module[\cC]{M}}{G_M}$ is a Quillen adjunction. The suspension and loop functors
        in $\cC$ induces a suspension functor $\Sigma_{\module{M}}$ and a loop functor $\Omega_{\module{M}}$ on $\module{M}$.
        We shall assume that, these induced suspension and loop functors commute with $G_M$.
    \end{enumerate}
\end{assumption} 

\begin{lemma}
    Under the assumptions \ref{asm:stabcomp}, 
    \[
      \adjunction{\Sigma_{M-mod}}{\module{M}}{\module{M}}{\Omega_{M-mod}}
    \]
    is a Quillen pair.
\end{lemma}
 \begin{proof}
     That $\Omega_{\module{M}}$ is a right Quillen functor is checked easily as follows. 
     We need to show that $\Omega_{\module{M}}$ preserves fibrations and acyclic fibrations.
     Let $f$ be a fibration (respectively, an acyclic fibration) in $\module{M}$. 
     Since both fibrations and weak equivalences in $\module{M}$ are created via $G_M$, 
     we need to check that $G_{M}(\Omega_{\module{M}}f)$ is a fibration (respectively, an acyclic fibration)
     in $\cC$. Since  $\Omega$ is a right Quillen functor on $\cC$ and since $G_{M}\Omega_{\module{M}}=\Omega G_{M}$, 
     $\Omega_{\module{M}} $ is a right Quillen functor. 
     
     Thus $(\Sigma_{\module{M}} \dashv \Omega_{\module{M}})$ forms a Qullen adjuntion (see, eg,
     \cite[lemma 1.3.4]{MR1650134}).
 \end{proof}

The following lemma describes a criterion for the Quillen adjuntion $(\Sigma_{\module{M}} \dashv \Omega_{\module{M}})$ to be a
Quillen equivalence.

\begin{lemma} \label{lem:shloopqe}
    Continuing with the same setup as in the previous lemma, if  $G_{M}$  preserves cofibrant objects 
    then the Quillen adjunction $\Sigma_{\module{M}} \dashv \Omega_{\module{M}}$ is a Quillen equivalence. 
\end{lemma}
\begin{proof}
  Note that $\Sigma_{\module{M}}$ is a Quillen equivalence \cite[definition 1.3.12]{MR1650134} if a map
  $\Sigma_{\module{M}}X \rightarrow Y$ is a weak equivalence in $\module{M}$ if and only if the corresponding map 
  $X \rightarrow \Omega_{\module{M}}Y$ is a weak equivalence in $\module{M}$ for all cofibrant objects $X$ and
  fibrant objects $Y$ in $\module{M}$.

  Recall that weak equivalences in $\module{M}$ with the transferred model structure are created by weak equivalences in $\cC$
  via $G_{M}$ and that $G_{M}\Sigma_{\module{M}}=\Sigma G_{M}$.  
  Let $x = G_M X$ and $y = G_M Y$. Since $G_{M}$ is a right Quillen functor, $y$ is a fibrant object in $\cC$,
  and by our hypothesis, $x$ is a cofibrant object in $\cC$.
  Therefore, we are reduced to showing that $\Sigma x \rightarrow y$ is a weak equivalence in $\cC$ if and only if
  the corresponding map $x \rightarrow \Omega y$ is a weak equivalence in $\cC$. This already holds as 
  $\Sigma \dashv \Omega$ is a Quillen equivalence on $\cC$.
\end{proof}

\begin{definition} \label{def:stabcomp}
  A monad $M$ satisfying the assumptions \ref{asm:stabcomp} is called \emph{stable compatible}
  (or, \emph{$M$ satisfies stable compatibility}) if the corresponding forgetful functor $G_M$
  preserves cofibrant objects.    
\end{definition}

We now highlight the importance of stable compatibility in the following theorem.
\begin{theorem} \label{the:stabcomp}
  Let $\cC$ be a stable cofibrantly generated model category.
  Suppose $M$ is a monad on $\cC$ which satisfies the following.
  \begin{enumerate}
    \item
      The category $\module{M}$ admits a right transferred model structure, which is 
      cofibrantly generated, and such that the induced adjunction $\Sigma_{\module{M}}
      \dashv \Omega_{\module{M}}$ is a Quillen adjunction. 
    \item The monad $M$ satisfies stable compatibility.
  \end{enumerate}
  Then $\module{M}$ is a stable model category. Moreover, if the monad is induced by a
  Quillen adjunction $\adjunction{F}{\cC}{\cD}{G}$, the underlying comparison functor $K$
  is such that $K \circ F=F_{M}$ and $G=G_{M} \circ K$ as in the diagram below   
  \[
    \begin{tikzcd}
      \cC  \arrow[ddrr,"F_{M}", shift left=1mm] \arrow[dd,"F", shift left=1mm] && 
      \\
      && \\
      \cD \arrow[dashrightarrow,rr,"\exists !K" below] \arrow[uu,"G",shift left=1mm] && \module{M} \arrow[uull, "G_{M}",shift left=1mm] 
    \end{tikzcd}
  \]
  and it preserves fibrations and acyclic fibrations.
\end{theorem}

\begin{proof}
  That $\module{M}$ is a stable model category is immediate from our previous discussion (see lemma \ref{lem:shloopqe}). 
   
   The existence of $K$ follows as usual. To see that $K$ preserves fibrations (respectively, acyclic fibrations),
   we use the fact that $G_{M}$ creates fibrations and weak equivalences in a right transferred model structure.
   Let $f$ be a fibration (respectively, acyclic fibration) in $\module{M}$. Since $G$ is a right Quillen functor, 
   $Gf$ is again a fibration (respectively acyclic fibration). Now, $G=G_{M} \circ K$ implies that $Kf$ is a 
   fibration (respectively, acyclic fibration).  
\end{proof}

\begin{remark} \label{rmk:stabcomp}
\begin{enumerate}
    \item The categorical diagram involving the Kleisli category as in Proposition \ref{pro:kleisliem} does not extend in this scenario, 
      as it is not a complete-cocomplete category.
    \item It is not true in general that $G_{M}$ commutes with suspension, i.e. $G_{M}\Sigma_{M-mod} \neq \Sigma G_{M}$ (however, see remark \ref{rem:gmstable} below).
    \item \label{itm:simmonad} It is a well-known result that if $\cC$ is a simplicial,  cofibrantly generated model category 
    with an action of a simplicial monad, which satisfies the assumptions of Schwede and Shipley \cite[lemma 2.3]{MR1734325},
    then $\module{M}$ is a simplicial cofibrantly generated model category with 
    $(F_{M} \dashv G_{M})$ being a simplicial Quillen adjunction (see Remark \ref{rem:postcofib}\eqref{rem:simppath}). The condition of stable compatibility
    extends this to a stable model structure.
    \item  \label{itm:chRasamc} 
    The converse to above theorem is not true in general; that is, the module category
    being stable does not imply that the monad is stable compatible. We give the
    following example.
    
    Consider the case of a commutative ring $R$ and an $R$-algebra $S$. We consider the projective/standard model structure on $Ch(R)$, the category of unbounded chain complexes of $R$-modules, as described in \cite[2.3.11]{MR1650134}. It is known to satisfy the monoid axiom of Schwede and Shipley (see \cite[section 5, \textit{Chain complexes}]{MR1734325}).
    Consider the monoid object given by the discrete DGA $S_{\bullet}= \cdots 0 \rightarrow S \rightarrow 0 \cdots$ in $Ch(R)$, where $S$ is situated at degree $0$.
    By definition, it follows that the Eilenberg-Moore category of 
    $S_{\bullet}$-modules in $Ch(R)$ is isomorphic to the category $Ch(S)$;
    the transferred model structure on the monad $T=(S_{\bullet} \bigotimes \_)$ 
    is the same as projective model structure on $Ch(S)$;
    and the forgetful functor $G_{M}$ is just degreewise restriction of scalars.
    Moreover, the suspension functor is the shift functor in this case
    and it commutes with forgetful functor. Further, it is a Quillen equivalence on both categories. By \cite[2.3.6]{MR1650134}, if $A_{\bullet}$ is a cofibrant object,
    each $A_{n}$ is a projective S-module for each $n$. However, restricting
    a projective $S$-module is not a projective $R$-module unless 
    $S$ is also a projective $R$-module. Thus $G_M$ does \emph{not} preserve cofibrant
    objects and stable compatibility fails.
\end{enumerate}
\end{remark}

\subsection{Some criteria for stable compatibility}
We would now like to point out some criteria under which the hypothesis of theorem \ref{the:stabcomp}. We'll use this to show that stable compatibility holds in the applications stated in this paper. The idea is to allow $G_{M}$ to have a good homotopical control. We state it in the form of the following lemma.
\begin{lemma} \label{lem:cndforsc}
    If the underlying endofunctor of the monad $M:\cC \rightarrow \cC$ preserves colimits, generating cofibrations, and generating acyclic cofibrations, then the monad satisfies hypothesis of theorem \ref{the:stabcomp}.
\end{lemma}
\begin{proof}
Let $I$ be the set of generating cofibrations, and $J$ the set of generating acyclic cofibrations of $\cC$. The candidate for
respective generating sets in a transferred model structure are exactly $F_{M} I$ and $F_{M} J$. Now, to transfer the right-induced
model structure, we check the conditions of Schwede and Shipley, as in lemma \ref{lem:schwship}.

Let $f: X \rightarrow Y$ be an $F_{M}(J)$-cell complex as in definition \ref{def:I-cell}, we need to show that $f$ is a weak equivalence. This is same as showing that
$G_{M} f$ is a weak equivalence. Let $f$ be the $\lambda$-transfinite composition of the following pushouts
\begin{equation} \label{eqn:pushoutf}
    \begin{tikzcd}
       \bigsqcup F_{M}A \arrow{r} \arrow[swap]{d} & X_{\beta} \arrow{d} \\
        \bigsqcup F_{M}B \arrow{r} & X_{\beta +1}
     \end{tikzcd}
\end{equation} 
where $X_{0}=X$.

Recall that $M=G_{M}F_{M}$ and $M$ preserves colimits (by assumption). Also $G_M$ creates, and hence preserves, colimits (cf \cite[theorem 5.6.12]{MR4727501} and \cite{MR2178101}).
Applying $G_{M}$ to the diagram \eqref{eqn:pushoutf}, we get the following $M(J)$-cell complexes in $\cC$
\begin{equation} \label{eqn:MJccmplx}
    \begin{tikzcd}
       \bigsqcup MA \arrow{r} \arrow[swap]{d} & G_{M}X_{\beta} \arrow{d} \\
        \bigsqcup MB \arrow{r} & G_{M}X_{\beta +1}
\end{tikzcd}
\end{equation}
Since $M$ preserves $J$, diagram \eqref{eqn:MJccmplx} is in fact a $J$-cell in $\cC$
and $G_{M} f$ is a $J$-cell in $\cC$ hence an acyclic cofibration and hence a weak equivalence.

Similarly, we can start with $X_{0}=*$ in $F_{M}I$-cell complexes and since $M$ preserves $I$ and $G_{M}$ preserves retract, we get the stable
compatibility of $M$.
\end{proof}

\begin{corollary} \label{cor:prescofb}
    The forgetful functor $G_{M}$ preserves cofibrations.
\end{corollary}
\begin{proof}
    Let $f:X \rightarrow Y$ be a cofibration in $\module{M}$. The proof of lemma \ref{lem:cndforsc} shows that $G_{M}(F_{M}I\text{-cell complex}) \subset I\text{-cell complex}$ in $\cC$. Since a cofibration in $\module{M}$ is a retract of an $F_{M}I$-cell complex and $G_{M}$ preserves both retracts and
    cell complexes, $G_{M}f$ is a retract of an $I$-cell complex in $\cC$. Hence, $G_M f$ is a cofibration.
\end{proof}
\begin{remark} \label{rem:gmstable}
    It is straightforward to check assumption \ref{asm:stabcomp}\eqref{itm:suscomGM} under hypothesis of \ref{lem:cndforsc}, 
    provided $\Sigma_{\module{M}}$ is well-defined and makes a Quillen pair with $(\Sigma_{\module{M}} \dashv \Omega_{\module{M}})$. 
    By definition, $\Sigma_{\module{M}}$X is the homotopy pushout of the diagram
    \[
    \begin{tikzcd}
      X \arrow{r} \arrow[swap]{d} & *  \\
       * 
    \end{tikzcd}
    \] 
   which is the ordinary pushout but with the above diagram replaced by a cofibrant one. 
   By corollary \ref{cor:prescofb}, $G_{M}$ preserves cofibrations, and we know that it preserves push out diagrams. Thus, it follows that
   $G_M$ commutes with suspension.
\end{remark}

\section{Separable monads} \label{sec:sepmonad}

In \cite{balmer:septrcat}, Balmer showed that for a separable monad on a triangulated category, the corresponding
category of modules $\module{M}$ admitted a triangulated structure, where the corresponding functor $G_M$ generated the
distinguished triangles. In this section, we study the effect of introducing a condition like separability in our model category set up.

We would like to recall the difference between the notion of pretriangulated in Hovey \cite{MR1650134} and Balmer \cite{balmer:septrcat}
to avoid confusion. For Hovey, a pretriangulated structure on a homotopy category of a model category $\cC$ means a collection of 
cofiber sequences and fiber sequences, which satisfy the
axioms as in definitions \ref{def:preTristruct} and \ref{def:pretriCat}. Such a structure always exists on the homotopy category of a 
cofibrantly generated model category (see remark \ref{def:Hoveytricat}). He then defines a triangulated category as a pretriangulated category when 
the suspension $\Sigma$ is an equivalence of homotopy category as in definition \ref{def:Hoveytricat}. Hovey observes that these triangulated categories are triangulated
in the classical sense described in the next paragraph.

For Balmer, a pretriangulated category means collection of distinguished triangles in a category satisfying all the axioms of
a triangulated category, except the octahedral axiom. See for example, \cite[Chapter 1, section 1.1]{neeman:trcat}. Here, the suspension functor
is always an autoequivalence (see definition \ref{def:monad}). A triangulated category \cite[definition 1.3.13 and proposition 1.4.6]{neeman:trcat} is a pretriangulated category where the octahedral axiom \cite[remark 1.4.7]{neeman:trcat} holds.

In this section, we study Balmer's treatment of separability for triangulated categories coming from the model categorical framework.
As mentioned in the introduction, this provides more information on \cite[remark 6.1]{balmer:septrcat}. In our situation, the distinguished triangles are the cofiber sequences of 
definitions \ref{def:cofseq} and \ref{def:boundarymap} as explained in remark \ref{rem:factstri}\eqref{ite:hoveyvsverdier}.
        
In the case of triangulated categories \cite{balmer:septrcat}, suspension functor on $\cC$ induce a suspension functor on
$\module[\cC]{M}$, which turns out to be autoequivalent. Contrast with this the situation involving a suitable model category $\cC$ with a monad $M$
satisfying stable compaibility. In this case, the induced suspension functor on $\Ho(\module[\cC]{M})$ is not automatically
an equivalence. It happens to be an autoequivalence, for example, when the suspension functor on the model category $\module[\cC]{M}$
is a Quillen equivalence. Hence, stable compatibility is much weaker than the assumptions in \cite{balmer:septrcat}.

Balmer further works with idempotent complete categories. We continue with the assumption that the categories
involved are idempotent complete in order to recover the work. 

\begin{definition}
    A model category is \emph{homotopy idempotent complete} if its homotopy category is an idempotent-complete category.
\end{definition}

\begin{definition}
    A monad $(M, \eta,\mu)$ is called \textit{separable} if the multiplication $\mu \colon M \circ M \rightarrow M$ admits section 
    $\sigma \colon M \rightarrow M \circ M$, that is $\mu \circ \sigma = \id_M$, and such that $M \mu \circ \sigma M = \sigma \circ \mu
    = \mu M \circ M \sigma$.
\end{definition}
\begin{remark}
  Note that, by Balmer \cite[section 3, \S\S 3.5 -- 3.11]{balmer:septrcat}, $M$ is separable, if and only if $G_M$ is a separable functor, which is equivalent to saying that the counit $\epsilon \colon F_{M}G_{M} \rightarrow \id_{\module{M}}$ admits a section $\upxi: \id_{\module{M}} \rightarrow F_{M}G_{M}$.
\end{remark}

\begin{assumption} \label{asm:sepassmpt}
For the rest of this section, we assume that the separable  monad also satisfies the assumptions of theorem \ref{the:stabcomp},
along with either the assumptions of lemma \ref{lem:cndforsc}, or remark \ref{rmk:stabcomp}\eqref{itm:simmonad} so that $G_M$ becomes a bi-Quillen functor.
\end{assumption}

\begin{example} \label{exa:sepalgeb}
    An example of such a monad comes alike from remark 
    \ref{rmk:stabcomp}\eqref{itm:chRasamc} when the $R$-algebra $S$ is finite \'etale . The forgetful functor $G_{M}$ becomes a 
    bi-Quillen functor and, in particular, satisfies stable compatibility. 
    The sections apply degreewise to get a 
    separable monad $S \bigotimes \_$ on $\Ch(R)$ with the standard/projective model structure mentioned in 
    \ref{rmk:stabcomp}\eqref{itm:chRasamc}.
\end{example}

We begin with a lemma that reflects the interaction of a separable monad with the transferred model structure.
We work with maps between bifibrant objects in $\module{M}$. For a map $f$ in the model category,, the corresponding map in the homotopy category
is denoted by $[f]$.  Using the small object arguement, factorize $f$ as $X \xrightarrow{i_{f}} E_{f} \xrightarrow{p_{f}}Y$, where $i_{f} \in F_{M}I\text{-cell}$ and $p_{f} \in RLP(F_{M}I)$, where $RLP(F_M I)$ denotes all the maps satisfying right lifing property with respect to maps in $F_M I$. 
As $X$ is cofibrant, the object $E_{f}$ is cofibrant  and as $Y$ is fibrant, $E_{f}$ is fibrant too. Therefore $E_{f}$ is bifibrant.

\begin{lemma} \label{lem:iftoFiGf}
    Under assumption \ref{asm:sepassmpt}, if $M$ is separable, then there exists a map in the arrow 
    category $i_{f} \rightarrow F_{M}i_{G_{M}f}'$ .    
\end{lemma}
\begin{proof}
 Let us consider the factorisation of $G_{M}f$ in $\cC$. Applying the small object argument in $\cC$, we get $G_{M}X \xrightarrow{i'_{G_{M}f}}E'_{G_{M}f} \xrightarrow{p'_{G_{M}f}}G_{M}Y$, along with the following commutative diagram in $\cC$:
\begin{equation} \label{eqn:cofibar1}
    \begin{tikzcd}
       G_{M}X \arrow{r}{i'_{G_{M}f}} \arrow{r} \arrow[swap]{d}{G_{M}i_{f}} & E'_{G_{M}f} \arrow{d}{p'_{G_{M}f}} \\
        G_{M}E_{f} \arrow{r}{G_{M}p_{f}} & G_{M}Y.
     \end{tikzcd}
\end{equation} 
By corollary \ref{cor:prescofb} or a consequence of $G_M$ being a bi-Quillen functor, the map $G_{M} i_{f}$ is a cofibration, and since the map $p'_{G_{M} f}$ is an acyclic fibration, by the RLP of acyclic fibrations with respect to cofibrations, we get a lift in the diagram \ref{eqn:cofibar1} which we denote by $t_{f}$:  
\begin{equation} \label{eqn:cofibar2}
    \begin{tikzcd}
       G_{M}X \arrow{r}{i'_{G_{M}f}} \arrow{r} \arrow[swap]{d}{G_{M}i_{f}} & E'_{G_{M}f} \arrow{d}{p'_{G_{M}f}} \\
        G_{M}E_{f} \arrow{r}{G_{M}p_{f}} \arrow[ur, dashed,"t_{f}"] & G_{M}Y
     \end{tikzcd}
\end{equation} 
 
Since both the maps $p'_{G_{M}f}$ and $G_{M}p_{f}$ are weak equivalences, by the 2-out-of-3 property of weak equivalences, so is $t_{f}$. Further, since $G_{M}X$ is a cofibrant object, so is $E'_{G_{M}f}$, therefore $t_{f}:G_{M}E_{f} \rightarrow E'_{G_{M}f}$ is a weak equivalence between cofibrant objects.

Since $F_{M}$ is a left Quillen functor, it preserves weak equivalences between cofibrant objects. Applying $F_{M}$ to the commutative diagram \eqref{eqn:cofibar2}, we get the following commutative diagram in $\module{M}$:
\begin{equation} \label{eqn:cofibar3}
    \begin{tikzcd}
       F_{M}G_{M}X \arrow{r}{F_{M}i'_{G_{M}f}} \arrow{r} \arrow[swap]{d}{F_{M}G_{M}i_{f}} & F_{M}E'_{G_{M}f} \arrow{d}{F_{M}p'_{G_{M}f}} \\
        F_{M}G_{M}E_{f} \arrow{r}{F_{M}G_{M}p_{f}} \arrow[ur, dashed,"F_{M}t_{f}"] & F_{M}G_{M}Y
     \end{tikzcd}
\end{equation} 

Now, by the upper-left triangle of diagram \eqref{eqn:cofibar3} and the naturality of sections
$\upxi_{X} \colon X \rightarrow F_{M}G_{M}X$ and $\upxi_{E_{f}} \colon E_{f} \rightarrow F_{M}G_{M}E_{f}$, 
the map $(\upxi_{X},(F_{M}t_{f}\circ \upxi_{E_{f}}))$ is the map claimed in the statement of the lemma:
\begin{equation}
    \begin{tikzcd}
       X \arrow[rr, "\upxi_{X}"]  \arrow[d, "i_{f}"'] & & F_{M}G_{M}X \arrow[d, "F_{M}i'_{G_{M}f}"] \\
       E_{f} \arrow[r, "\upxi_{E_{f}}"] 
       \arrow[rr, bend right, "F_M t_f \circ \upxi_{E_f}"'] & F_M G_M E_f \arrow[r, "F_{M}t_{f}"] & F_{M} E'_{G_{M}f}
     \end{tikzcd}
\end{equation} 
\end{proof}

\begin{remark}
We note that $X$, $E_{f}$, $F_{M}G_{M}X$, $F_{M}E'_{G_{M}f}$ are all cofibrant objects in $\module{M}$
and $i_{f}$ and $F_{M}i'_{G_{M}f}$ are cofibrations.
\end{remark}

Next, we want to extend the above map into a map of cofiber sequences. Before that, we recall a result from Hovey.
\begin{lemma}[Hovey \cite{MR1650134}, Proposition 6.2.5] \label{lem:extcofib}
Suppose $\cC$ is a pointed model category and we have a commutative square of cofibrant objects
\[
  \begin{tikzcd}
        A' \arrow{r}{f'} \arrow{d}{q_{1}} \arrow[swap]{d} & B' \arrow{d}{q_{2}} \\
        A \arrow{r}{f} & B
  \end{tikzcd}
\] 
where $f'$ and $f$ are cofibrations, with cofibers $g':B' \rightarrow C'$ and $g: B \rightarrow C$ respectively. 
Then the induced map $q_{3}:C' \rightarrow C$ is equivariant in $Ho(\cC)$ with respect to the cogroup homomorphism $\Sigma q_{1}$.

The corresponding statement for fibrations hold by duality.   
\end{lemma}
\begin{corollary} \label{cor:indcisom}
    If $q_{1}$ and $q_{2}$ are weak equivalences in lemma \ref{lem:extcofib} then the induced map between the
    cofiber sequences is an isomorphism.
\end{corollary}
\begin{proof}
Since the cofiber of  $f': A' \rightarrow B'$ map is the pushout of the diagram: 
\begin{equation} \label{eqn:diaapcof}
    \begin{tikzcd}
       A' \arrow{r}{f'}  \arrow[swap]{d} & B'  \\
        * 
     \end{tikzcd}
\end{equation} 

and similary the cofiber of $f: A \rightarrow B$:
\begin{equation} \label{eqn:diagacof}
    \begin{tikzcd}
       A \arrow{r}{f}  \arrow[swap]{d} & B \\
        * 
     \end{tikzcd}
\end{equation} 
    Then $(\id_{*},q_{1},q_{2})$ are the object-wise weak equivalences between the diagrams \eqref{eqn:diaapcof} and 
    \eqref{eqn:diagacof}. By the cube lemma \cite[lemma 5.2.6]{MR1650134}, the induced map $q_{3}: C' \rightarrow C$
    is a weak equivalence as well.
\end{proof}

\begin{lemma} \label{lem:dirincof}
The map $i_{f} \rightarrow F_{M}i_{G_{M}f}'$ in lemma \ref{lem:iftoFiGf} extends to an inclusion map
of the cofiber sequence $C_1$ associated with $i_{f}$ into the cofiber sequence $C_2$ associated with
$F_{M}i_{G_{M}f}'$ in $\Ho(\module{M})$. Moreover, the sequence $C_1$ is a direct summand of $C_2$ under this inclusion map.
\end{lemma}
\begin{remark}
    Before we proceed to give a proof, note that this lemma is not saying that cofiber sequence associated with $i_{f}$ in $\Ho(\module{M})$is a direct summand of $L F_{M}$ applied to the cofiber sequence of $i'G_{M}f$ in $\Ho \cC$. Nor is it implying that the triangle associated to $[f]$ is a direct summand of the $LF_{M}$ applied to the triangle associated to the map $RG_{M}([f])$. For the latter follow Balmer's proof of 
    separability \cite{balmer:septrcat}.
\end{remark}
\begin{proof}
We note that the map $i_{f} \rightarrow F_{M}i_{G_{M}f}'$ as defined factors as 
    \[\begin{tikzcd}
         X \arrow{r}{\upxi_{X}} \arrow[swap]{d}{i_{f}} & F_{M}G_{M}X \arrow{d}{F_{M}G_{M}i_{f}} \arrow{r}{id} & F_{M}G_{M}X \arrow{d}{F_{M}i'_{G_{M}f}}  \\
         E_{f} \arrow{r}{\upxi_{E_{f}}} & F_{M}G_{M}E_{f} \arrow{r} {(i_{2}=F_{M}t_{f})} & F_{M}E'_{G_{M}f} 
    \end{tikzcd}\]

By naturality of sections and counit, the first square \[\begin{tikzcd}
         X \arrow{r}{\upxi_{X}} \arrow[swap]{d}{i_{f}} & F_{M}G_{M}X \arrow{d}{F_{M}G_{M}i_{f}} \arrow[l,"\epsilon_{X}",shift left=1ex]   \\
         E_{f} \arrow{r}{\upxi_{E_{f}}} & F_{M}G_{M}E_{f}  \arrow[l,"\epsilon_{E_{f}}",shift left=1ex]
    \end{tikzcd}\]   

gives an idempotent of the square
\[\begin{tikzcd}
         F_{M}G_{M}X \arrow{r}{\upxi_{X}\epsilon_{X}} \arrow[swap]{d}{F_{M}G_{M}i_{f}} & F_{M}G_{M}X \arrow{d}{F_{M}G_{M}i_{f}}    \\
         F_{M}G_{M}E_{f} \arrow{r}{\upxi_{E_{f}}\epsilon_{E_{f}}} & F_{M}G_{M}E_{f}  
    \end{tikzcd}\] 
by lemma \ref{lem:extcofib}, this extends to an idempotent map of the cofiber sequence associated to the map $F_{M}G_{M}i_{f}$ and as $\Ho(\module{M})$ is an 
idempotent complete category, the cofiber sequence associated to $i_{f}$ is a direct summand of the former. Moreover, the pair 
$(\id_{F_{M}G_{M}X}, F_{M}t_{f})$ is a pair of weak equivalences, and by corollary \ref{cor:indcisom}, the cofiber sequence associated to $F_{M}G_{M}i_{f}$ is 
isomorphic to the cofiber sequence associated with the map $F_{M}i'_{G_{M}f}$ . So via the isomorphism the cofiber associated to $i_{f}$ becomes a direct summand of 
the cofiber associated to $F_{M}i'_{G_{M}f}$.
\end{proof}

\begin{remark} \label{rem:sepcompu}
Let us consider a pair of composable maps  $X \xrightarrow{f} Y \xrightarrow{g} Z$. Applying the functorial factorization, we obtain 
$i_{f}$, $i_{g}$ and $i_{gf}$ and since the factorization is functorial on the arrow category, the map $(\id_{X},g)$  
\[\begin{tikzcd}
        X \arrow{r}{f} \arrow{d}{\id_{X}}  & Y \arrow{d}{g} \\
        X \arrow{r}{gf} & Z
     \end{tikzcd}
     \]   
induces a map $s_{f}:E_{f} \rightarrow E_{gf}$ 
\[\begin{tikzcd}
        X \arrow{r}{i_{f}} \arrow{d}{\id_X}  & E_{f} \arrow{d}{s_{f}} \arrow{r}{p_{f}} & Y \arrow{d}{g}\\
        X \arrow{r}{i_{gf}} & E_{gf} \arrow{r}{p_{gf}} & Z
     \end{tikzcd}
     \]   
such that $s_{f}i_{f}=i_{gf}$ and $gp_{f}=p_{gf}s_{f}$ and similarly the map $(f,\id_{Z})$  
\[\begin{tikzcd}
        X \arrow{r}{gf} \arrow{d}{f}  & Z \arrow{d}{id_{Z}} \\
        Y \arrow{r}{g} & Z
     \end{tikzcd}
     \]  
     induces a map $s_{gf}:E_{gf} \rightarrow E_{g}$ giving a commutative diagram  
     \[\begin{tikzcd}
        X \arrow{r}{i_{gf}} \arrow{d}{f}  & E_{gf} \arrow{d}{s_{gf}} \arrow{r}{p_{gf}} & Z \arrow{d}{\id_{Z}}\\
        Y \arrow{r}{i_{g}} & E_{g} \arrow{r}{p_{g}} & Z.
     \end{tikzcd}
     \]
     Note that we obtain similar commutative diagrams in $\cC$ when we apply the functorial factorization to $G_{M}f$, 
     $G_{M}g$ and $G_{M}gf$. We denote the middle maps by $s'_{G_{M}f}:G_{M}X \rightarrow E'_{G_{M}f}$ 
     and $s'_{G_{M}gf}:G_{M}X \rightarrow E'_{G_{M}gf}$ .
     
    Now the whole diagram 
    \begin{equation} \label{eqn:octahed1}
       \begin{tikzcd}
        X \arrow{r}{i_{f}} \arrow{d}{\id_{X}}  & E_{f} \arrow{d}{s_{f}}  \\
        X \arrow{r}{i_{gf}} \arrow{d}{f} & E_{gf} \arrow{d}{s_{gf}} \\
        Y \arrow{r}{i_{g}} & E_{g}
       \end{tikzcd}
    \end{equation} 
      admits sections into 
      \begin{equation} \label{eqn:octahed2}
        \begin{tikzcd}
        F_{M}G_{M}X \arrow{r}{F_{M}G_{M}i_{f}} \arrow{d}{\id_{F_{M}G_{M}X}}  & F_{M}G_{M}E_{f} \arrow{d}{F_{M}G_{M}s_{f}}  \\
        F_{M}G_{M}X \arrow{r}{F_{M}G_{M}i_{gf}} \arrow{d}{F_{M}G_{M}f} & F_{M}G_{M}E_{gf} \arrow{d}{F_{M}G_{M}s_{gf}} \\
        F_{M}G_{M}Y \arrow{r}{F_{M}G_{M}{i_{g}}} & F_{M}G_{M}E_{g}
     \end{tikzcd}
      \end{equation} 
      by naturality of the sections, thereby extending into a direct summand of ocathedron associated to \eqref{eqn:octahed1} into oactahedron associated to 
      \eqref{eqn:octahed2}. Furthermore since the pairs $(\id_{F_{M}G_{M}X}, F_{M}t_{f})$, $(\id_{F_{M}G_{M}X}, F_{M}t_{gf})$ 
      and $(\id_{F_{M}G_{M}Y},  F_{M}t_{g})$ all satisfy hypothesis of corollary \ref{cor:indcisom} the octahedron associated to \eqref{eqn:octahed2} is isomorphic to the octahedron associated to \eqref{eqn:octahed3}. So the octahedron associated to \eqref{eqn:octahed1} becomes a direct summand of octahedron associated to \eqref{eqn:octahed3}.
    \begin{equation} \label{eqn:octahed3}
     \begin{tikzcd}
        F_{M}G_{M}X \arrow{r}{F_{M}i'_{G_{M}f}} \arrow{d}{\id_{F_{M}G_{M}X}}  & F_{M}E'_{G_{M}f} \arrow{d}{F_{M}s'_{G_{M}f}}  \\
        F_{M}G_{M}X \arrow{r}{F_{M}i'_{G_{M}gf}} \arrow{d}{F_{M}G_{M}f} & F_{M}E'_{G_{M}gf} \arrow{d}{F_{M}s'_{G_{M}gf}} \\
        F_{M}G_{M}Y \arrow{r}{F_{M}i'_{G_{M}g}} & F_{M}E'_{G_{M}g}
     \end{tikzcd}
    \end{equation} .
\end{remark}

We summarize the computation done in the remark above as the following theorem.

\begin{theorem} \label{the:mainsepm}
Let $\cC$ be a cofibrantly generated stable model category, and let $M$ be a 
separable monad on $\cC$ such that $M$ satisfies the assumptions of theorem 
\ref{the:stabcomp} along with the assumption \ref{asm:sepassmpt}. Then $\module{M}$ is a cofibrantly generated stable model 
category. Moreover, a cofiber sequence admits a direct summand inclusion into a suitable
cofiber sequence.
\end{theorem}
\begin{proof} For the moreover part, we note that if $[f]:X  \rightarrow Y$ is a map in $\Ho(\module{M})$. We can represent it by some $\Tilde{f}: \Tilde{X} \rightarrow \Tilde{Y}$ in $\module{M}$, where $\Tilde{X}$ and $\Tilde{Y}$ are bifibrant objects in $\module{M}$. We conclude by lemma \ref{lem:dirincof} and remark \ref{rem:sepcompu}.
\end{proof}
\begin{remark} \label{rem:sepscmnd}
    \begin{enumerate}
        \item With the availability of sections along with stable compatibility, we note that to show a cofiber sequence  admits a direct summand inclusion as above, we did not appeal to the trick of lifting idempotents up to nilpotence as in \cite[lemma 1.6(c)]{balmer:septrcat}; yet we obtain a similar looking conclusion about the data of cofiber sequences-distinguished triangles as in \cite[theorem 4.1]{balmer:septrcat}.
        \item  By induction, we can now extend the same observation as above to an $n$-triangle (see \cite[definition 5.1]{balmer:septrcat}) for $n \ge 1$ as noted in \cite[remark 5.14]{balmer:septrcat} and again obtain similar looking conclusion as in \cite[theorem 5.17]{balmer:septrcat}.
    \end{enumerate}
\end{remark}
\section{Applications} \label{sec:applications}
For all the applications in this paper we are in the situation where a functor $G: \cD \rightarrow \cC$ admits both left and right adjoints i.e. $(F \dashv G \dashv H)$, this gives us two monads: $M=GF \colon \cC \rightarrow \cC$ and  $N=HG \colon \cD \rightarrow \cD$. Let $\module{M}$ and $\module{N}$ be the respective module categories. By the adjoint lifting theorem, we get $(F' \dashv G' \dashv H')$ the induced adjoints for $G': \module{N} \rightarrow \module{M}$. 

\subsection{Category of modules over related monads} \label{ssc:drumhack}

Here, we recall and improve upon a result of Drummond-Cole and Hackney.  We describe the similarities of properties shared by $\module{M}$ and $\module{N}$.

    We recall (Theorem 2.3 \cite{MR3937590}): Suppose $(\mathcal{M},\mathcal{W},I,J)$ is a cofibrantly generated model category and $\mathcal{N}$ is a bicomplete category. Let $F: \mathcal{N} \rightarrow \mathcal{M}$ be a functor which is right adjoint to $L$ and left adjoint to $R$. If $(FL,FR)$ is a Quillen adjunction, then there exists a right-induced cofibrantly generated model category structure $(\mathcal{N},F^{-1}\mathcal{W},LI,LJ)$. Both $(L,F)$ and $(F,R)$ are Quillen adjunctions.

    \begin{remark} \label{rem:autostcy}
    We point out that stable compatibility is always satisfied when the functor that creates the right induced structure (for example, for $G_{M}$ in section \ref{sec:stablecpt} in page \pageref{itm:rtmodstr}, and $F$ in theorem stated above), and is both a left and a right Quillen functor. So if $\mathcal{M}$ is stable then $\mathcal{N}$ is stable as well.
    \end{remark}
    
Let  $A$ be a cofibrant highly structured ring spectrum ($S$-module)
as defined in \cite{MR1417719} and $f \colon S \rightarrow A$ be the unit map. 
Denote by 
$\mathcal{M}_{S}$ the cofibrantly generated topological (and hence simplicial), stable model category of $S$-modules and 
let $\mathcal{M}_{A}$ be the cofibrantly generated topological (and hence 
simplicial), stable model category of $A$-modules.  The unit map induces 
$(f_{*}=A \wedge_{S}-) \dashv f^{*} \dashv F_{S}(A,-))$.

The functor $f^{*}$ is both a left and right Quillen functor. For $M=f^{*}f_{*}$, $
\module{M} = \mathcal{M}_{A}$ and so $F_{M}=f_{*}$ and $G_{M}=f^{*}$. We point out 
here that if $G_{M}$ is both a left and right Quillen functor, then it automatically 
satisfies our condition of stable compatibility.

So $\module{M}$ is always a cofibrantly generated simplicial and stable model category for any $A$. But, can we say anything about the category $\module{N}$?
So, if $(F' \dashv G' \dashv H')$ induced by $(f_{*}=A \wedge_{S}-) \dashv f^{*} \dashv F_{S}(A,-))$  is such that the triple satisfies \cite[theorem 2.3]{MR3937590}; and hence by remark \ref{rem:autostcy}, $\module{N}$ is also a cofibrantly generated simplicial and stable model category.

\subsection{Finite, flat maps between schemes} \label{ssc:finflmor}
 We now consider the monads in Algebraic geometry. We state the example for Noetherian schemes.
 
 Let $f: X \rightarrow Y$ be a finite, flat morphism of Noetherian schemes, for such a suitable $f$ we get an adjunction 
 $\adjunction{f^{*}}{\Ch(QCoh(Y))}{Ch(QCoh(X))}{f_{*}}$ and the monad $M_{f}=f_{*}f^{*}$ on $Ch(QCoh(Y))$.
We consider the locally-free model structure on $\chA$ as stated in theorem \ref{thm:locallyfreeModelStruct} for $\mathscr{A} = \QCoh(Y)$.
We claim that when $f$ is finite and flat, it satisfies the assumptions of lemma \ref{lem:cndforsc}. To see this, we note that we have the adjunctions $(f^{*} \dashv f_{*} \dashv f^{!} )$, and therefore $M_{f}=f_{*}f^{*}$ preserves colimits. The pull back $f^*$ is known to preserve locally free sheaves, and so does $f_*$ for finite (hence affine) and flat $f$ (see, for example, 
\cite[\href{https://stacks.math.columbia.edu/tag/02K9}{Tag 02K9}, lemma 29.49.2]{stacks-project}; and also \cite[II, exercise 5.8; III, proposition 9.2]{hartshorne:alggeom} for the case of reduced, connected noetherian schemes).
Hence, $M_f$ takes locally free sheaves to locally free sheaves. The functor $M_f$ also preserve
subobjects as both $f^*$ and $f_*$ are exact under the assumptions.
Therefore, $M_f$ preserves the set of generating cofibrations and generating acyclic 
cofibrations as mentioned in \ref{def: I and J}. Further, by \cite[exercise III, 8.2]{hartshorne:alggeom}, we know that 
for affine $f$, $R^i f_* = 0$ for $i > 0$, and hence $R f_* = f_*$

Thus, the monad $M_f = f_* f^*$ on $\chA$ is
computed as $M_f(A_{\bullet}) = (\dotsc \rightarrow f_* f^*(A_i) \rightarrow f_* f^* 
(A_{i+1}) \rightarrow \dotsc)$. From this, it is clear that elements of $\module[\chA]
{M_f}$ pairs $(A_{\bullet}, \lambda_{\bullet})$ with $\lambda_n \colon f_* f^* A_n 
\longrightarrow A_n$, where all the squares in the corresponding diagram commute.
The above description tells us that $\module[\chA]{M_f}$ is the same as 
$\Ch(\module[\cA]{M_f})$.  
Clearly, in this case, $T$ is the shift functor and hence it commutes with the
forgetful functor $G_M$. 
By theorem \ref{the:stabcomp}, we get that $\module[\chA]{M_f}$ is a 
cofibrantly generated stable model category.

Examples of such maps are \'etale maps \cite{balmer:etalneth} and the 
$n$-power map $z \mapsto z^n$.

\subsection{A remark about comonads and comodules over them} \label{ssc:comodules}

By considering monads in the opposite category, one gets the notion of 
comonads, and one can construct comodules over them such that the opposite 
category of comodules over a comonad $W$ over $\cC$ coincide with the module 
category of $W^{\text{op}}$ over $\cC^{\text{op}}$ 
(compare \cite[remark 5.19]{balmer:septrcat}). 
Thus, by considering the corresponding notion of stable compatibility, namely
that in the situation of a left-transferred model structure the suspension should commute with the forgetful functor from the 
category of comodules, and that the forgetful functor should preserve fibrant
objects, we would get a criterion for the category of comodules to be stable.

\bibliography{references}

\providecommand{\bysame}{\leavevmode\hbox to3em{\hrulefill}\thinspace}
\providecommand{\MR}{\relax\ifhmode\unskip\space\fi MR }
\providecommand{\MRhref}[2]{%
  \href{http://www.ams.org/mathscinet-getitem?mr=#1}{#2}
}
\providecommand{\href}[2]{#2}
\begin{thebibliography}{EKMM97}

\bibitem[Bal11]{balmer:septrcat}
Paul Balmer, \emph{Separability and triangulated categories}, Adv. Math.
  \textbf{226} (2011), no.~5, 4352--4372. \MR{2770453 (2012c:18021)}

\bibitem[Bal14]{balmer:etalneth}
Paul Balmer, \emph{The derived category of an étale extension and the
  separable neeman-thomason theorem}, 2014.

\bibitem[BW05]{MR2178101}
Michael Barr and Charles Wells, \emph{Toposes, triples and theories}, Repr.
  Theory Appl. Categ. (2005), no.~12, x+288, Corrected reprint of the 1985
  original [MR0771116]. \MR{2178101}

\bibitem[DCH19]{MR3937590}
Gabriel~C. Drummond-Cole and Philip Hackney, \emph{A criterion for existence of
  right-induced model structures}, Bull. Lond. Math. Soc. \textbf{51} (2019),
  no.~2, 309--326. \MR{3937590}

\bibitem[DS18]{MR3851535}
Ivo Dell'Ambrogio and Beren Sanders, \emph{A note on triangulated monads and
  categories of module spectra}, C. R. Math. Acad. Sci. Paris \textbf{356}
  (2018), no.~8, 839--842. \MR{3851535}

\bibitem[EKMM97]{MR1417719}
A.~D. Elmendorf, I.~Kriz, M.~A. Mandell, and J.~P. May, \emph{Rings, modules,
  and algebras in stable homotopy theory}, Mathematical Surveys and Monographs,
  vol.~47, American Mathematical Society, Providence, RI, 1997, With an
  appendix by M. Cole. \MR{1417719}

\bibitem[Har77]{hartshorne:alggeom}
Robin Hartshorne, \emph{Algebraic geometry}, Springer, 1977.

\bibitem[Hov99]{MR1650134}
Mark Hovey, \emph{Model categories}, Mathematical Surveys and Monographs,
  vol.~63, American Mathematical Society, Providence, RI, 1999. \MR{1650134}

\bibitem[Hov01]{MR1814077}
\bysame, \emph{Model category structures on chain complexes of sheaves}, Trans.
  Amer. Math. Soc. \textbf{353} (2001), no.~6, 2441--2457. \MR{1814077}

\bibitem[ML98]{maclane:catworkmath}
Saunders Mac~Lane, \emph{Categories for the working mathematician}, second ed.,
  Graduate Texts in Mathematics, vol.~5, Springer-Verlag, New York, 1998.
  \MR{1712872 (2001j:18001)}

\bibitem[Nee01]{neeman:trcat}
Amnon Neeman, \emph{Triangulated categories}, Annals of Mathematics Studies,
  vol. 148, Princeton University Press, Princeton, NJ, 2001. \MR{1812507}

\bibitem[Rie16]{MR4727501}
Emily Riehl, \emph{Category theory in context}, Aurora Dover Modern Math
  Originals, Dover Publications, Inc., Mineola, NY, 2016. \MR{4727501}

\bibitem[RV15]{MR3338539}
Emily Riehl and Dominic Verity, \emph{Completeness results for quasi-categories
  of algebras, homotopy limits, and related general constructions}, Homology
  Homotopy Appl. \textbf{17} (2015), no.~1, 1--33. \MR{3338539}

\bibitem[SS00]{MR1734325}
Stefan Schwede and Brooke~E. Shipley, \emph{Algebras and modules in monoidal
  model categories}, Proc. London Math. Soc. (3) \textbf{80} (2000), no.~2,
  491--511. \MR{1734325}

\bibitem[{Sta}18]{stacks-project}
The {Stacks Project Authors}, \emph{\textit{Stacks Project}},
  \url{https://stacks.math.columbia.edu}, 2018.

\end{thebibliography}
\bibliographystyle{amsalpha}

\end{document}